\documentclass[12pt]{article}
\def\q{\hfill\rule{1ex}{1ex}}
\def\0{\emptyset}

\def\p{{\bf Proof.} \quad}
\def\q{\hfill\rule{1ex}{1ex}}

\def\n{\noindent}

\usepackage{graphicx}

\begin{document}
\title{\bf Edge-fault-tolerant strong Menger edge connectivity of bubble-sort star graphs}
\author{{\small\bf
Jia Guo$^{1,}$$^2$}\thanks{email: guojia199011@163.com}\\
%\quad\quad
%{\small\bf Mei Lu$^3$}\thanks{email: lumei@mail.tsinghua.edu.cn}\\
{\small $^1$School of Software, Northwestern Polytechnical University, Xi'an, Shaanxi 710072, PR China}\\
{\small $^2$College of Science, Northwest A$\&$F University, Yangling, Shaanxi 712100, PR China}\\
%{\small $^3$Department of Mathematical Sciences, Tsinghua University, Beijing 100084, PR China}
}

\date{}
\maketitle\baselineskip 15.5pt

\begin{abstract}
\baselineskip=0.5cm
The connectivity and edge connectivity of interconnection network determine the fault tolerance of the network. An interconnection network is usually viewed as a connected graph, where vertex corresponds processor and edge corresponds link between two distinct processors. Given a connected graph $G$ with vertex set $V(G)$ and edge set $E(G)$, if for any two distinct vertices $u,v\in V(G)$, there exist $\min\{d_G(u),d_G(v)\}$ edge-disjoint paths between $u$ and $v$, then $G$ is strongly Menger edge connected. Let $m$ be an integer with $m\geq1$. If $G-F_e$ remains strongly Menger edge connected for any $F_e\subseteq E(G)$ with $|F_e|\leq m$, then $G$ is $m$-edge-fault-tolerant strongly Menger edge connected. If $G-F_e$ is strongly Menger edge connected for any $F_e\subseteq E(G)$ with $|F_e|\leq m$ and $\delta(G-F_e)\geq2$, then $G$ is $m$-conditional edge-fault-tolerant strongly Menger edge connected. In this paper, we consider the $n$-dimensional bubble-sort star graph $BS_n$. We show that $BS_n$ is $(2n-5)$-edge-fault-tolerant strongly Menger edge connected for $n\geq3$ and $(6n-17)$-conditional edge-fault-tolerant strongly Menger edge connected for $n\geq4$. Moreover, we give some examples to show that our results are optimal.

\vskip 0.3cm

%{\bf AMS} classification: 05C50  \vskip 0.3cm

{\bf Keywords:} fault-tolerance, strong Menger edge connectivity, bubble-sort star graph

\end{abstract}

\vskip.3cm

\n{\large\bf 1.\quad Introduction}
\vskip.2cm

The connectivity and edge connectivity are two crucial factors for the interconnection
networks since they determine the fault tolerance of the networks.
An interconnection network can be viewed as a simple connected graph, where vertex corresponds processor and edge corresponds link. In the rest of this paper, we only consider simple connected graphs and we follow the work of
\cite{Bondy1} for definitions and notations not defined here.

Let $G=(V(G),E(G))$ be a simple connected graph. For a vertex $v\in V(G)$, $N_G(v)=\{u~|~(u,v)\in E(G)\}$ is the set of neighbours of $v$ and $E_G(v)=\{(u,v)~|~(u,v)\in E(G)\}$ is the set of edges that are incident with $v$.
Let $d_G(v)=|N_G(v)|$ be
the {\em degree} of $v$ and $\delta(G)=\min\{d_G(v)~|~v\in V(G)\}$ be the {\em minimum degree} of $G$.
If $d_G(v)=k$ for every $v\in V(G)$, then $G$ is {\em $k$-regular}. $G$ is bipartite if there exist two vertex subsets $V_1,~V_2$ with $V_1\cap V_2=\emptyset$ such that $V(G)=V_1\cup V_2$ and for each edge $(u,v)\in E(G)$, $|\{u,v\}\cap V_1|=|\{u,v\}\cap V_2|=1$. It is well known that bipartite graphs contain no odd cycles.
Let $F_1,F_2\subseteq V(G)$ with $F_1\cap F_2=\emptyset$, denote $E_G(F_1, F_2)=\{(u,v)\in E(G)~|~u\in F_1,~v\in F_2\}$.
Let $F\subseteq V(G)$ and $F_e\subseteq E(G)$. We use $G-F$ to denote the subgraph of $G$ with vertex set $V(G)-F$ and edge set $E(G)-\{(u,v)\in E(G)~|~\{u,v\}\cap F\neq \emptyset \}$. If $G-F$ is disconnected or has only one vertex, then $F$ is a {\em vertex cut} of $G$. We use $G-F_e$ to denote the subgraph of $G$ with vertex set $V(G)$ and edge set $E(G)-F_e$. If $G-F_e$ is disconnected, then $F_e$ is an {\em edge cut} of $G$.
The {\em connectivity} (resp. {\em edge connectivity}) of $G$, denoted by $\kappa(G)$ (resp. $\lambda(G)$), is the minimum size of $F$ (resp. $F_e$) such that $F$ (resp. $F_e$) is a vertex cut (resp. an edge cut) of $G$.
$P_k=uv_2v_3\cdots v_{k-1}v$
on $k$ distinct vertices $u,v_2,\cdots,v_{k-1},v$ of $G$
is a {\em $(u,v)$-path} if $(u,v_2)\in E(G)$, $(v_{k-1},v)\in E(G)$, and $(v_i,v_{i+1})\in E(G)$ for every $i\in\{2,\cdots,k-2\}$.
$F\subseteq V(G)-\{u,v\}$ (resp. $F_e \subseteq E(G)$) is an {\em $(u,v)$-cut} (resp. {\em $(u,v)$-edge-cut}) if $G-F$ (resp. $G-F_e$) has no $(u,v)$-path. Menger's theorem is a classical theorem about the connectivity and edge connectivity.

\vskip.2cm

{\bf Theorem 1.1 \cite{Menger}} {\em Let $G$ be a graph and $u,v\in V(G)$ with $u\neq v$. Then

(1) the minimum size of an $(u,v)$-cut equals to the maximum number of disjoint $(u,v)$-paths for $(u,v)\not\in E(G)$;

(2) the minimum size of an $(u,v)$-edge-cut equals to the maximum number of edge-disjoint $(u,v)$-paths.}

\vskip.2cm

Motivated by Menger's theorem, Oh et al. \cite{Oh} proposed the strong Menger connectivity (also called the maximal local-connectivity) and Qiao et al. \cite{Qiao} introduced the strong Menger edge connectivity, which are showed in the following definition.
\vskip.2cm

{\bf Definition 1.2} {\em Let $G$ be a connected graph and $u,v\in V(G)$ be any two distinct vertices. Then

(1) $G$ is strongly Menger connected if there exist $\min\{d_G(u),d_G(v)\}$ disjoint $(u,v)$-paths;

(2) $G$ is strongly Menger edge connected if there exist $\min\{d_G(u),d_G(v)\}$ edge-disjoint $(u,v)$-paths.}

\vskip.2cm

Since edge faults may occur in real interconnection networks, the edge-fault-tolerant strong Menger edge connectivity has been proposed.
\vskip.2cm

{\bf Definition 1.3} {\em Let $m\geq1$ be an integer, $G$ be a connected graph, and $F_e\subseteq E(G)$ be any arbitrary edge subset of $G$ with $|F_e|\leq m$. Then

(1) $G$ is $m$-edge-fault-tolerant strongly Menger edge connected if $G-F_e$ is strongly Menger edge connected;

(2) $G$ is $m$-conditional edge-fault-tolerant strongly Menger edge connected if $G-F_e$ is strongly Menger edge connected for any $F_e$ with $\delta(G-F_e)\geq2$.}

\vskip.2cm

The edge-fault-tolerant strong Menger edge connectivity of many interconnection networks has been studied. For example, Qiao et al. proved that the folded hypercube is $(2n-2)$-conditional edge-fault-tolerant strongly Menger edge connected \cite{Qiao}. Li et al. discussed the edge-fault-tolerant strong Menger edge
connectivity of the hypercube-like network \cite{P.Li} and the balanced hypercube \cite{P.Li1}. He et al. considered the strong Menger edge connectivity of the regular network \cite{He}.

This paper deals with the edge-fault-tolerant strong Menger edge connectivity of the $n$-dimensional bubble-sort star graph $BS_n$ \cite{Chou}, which gains many nice properties, such as vertex transitive and high degree of
regularity.
Cai et al. showed that $BS_n$ is $(2n-5)$-fault-tolerant strongly Menger connected \cite{Cai}.
Wang et al. studied the 2-extra diagnosability \cite{Wang1}, the 2-good-neighbor diagnosability \cite{Wang2}, and the strong connectivity \cite{Wang3} of $BS_n$. Gu et al. discussed the pessimistic diagnosability of $BS_n$ \cite{Gu}. Zhao et al. investigated the generalized connectivity of $BS_n$ \cite{Zhao}.
Zhu et al. gave an algorithm to determine the $h$-extra connectivity of $BS_n$ of low dimensions \cite{Zhu}. Zhang et al. considered the structure connectivity and substructure connectivity of $BS_n$ \cite{Zhang}.

The remainder of this paper is organized as follows: Section 2 introduces the definition of $BS_n$ and gives some properties of $BS_n$. In section 3, we demonstrate the edge-fault-tolerant strong Menger edge connectivity of $BS_n$. In section 4, we discuss the conditional edge-fault-tolerant strong Menger edge connectivity of $BS_n$. Section 5 concludes this paper.

\vskip.3cm

\n{\large\bf 2.\quad Preliminaries}

\vskip.2cm

Let $l_1,l_2$ be two integers with $1\leq l_1\leq l_2$. Set $[l_1,l_2]=\{l~|~l_1\leq l \leq l_2, ~l~ \mbox{is an integer}\}$. Now we give the definition of the $n$-dimensional bubble-sort star graph $BS_n$.

\vskip.2cm
{\bf Definition 2.1 \cite{Chou}} The $n$-dimensional bubble-sort star graph $BS_n$ has
vertex set $V(BS_n)$ and edge set $E(BS_n)$. A vertex $v\in V(BS_n)$ if and only if $v$ is a permutation on $[1,n]$, which is
denoted as $v=
v_1v_2\cdots v_n$. Let $x=x_1x_2\cdots x_n\in V(BS_n)$, $y=y_1
y_2\cdots y_n\in V(BS_n)$ with $x\neq y$. Then $(x,y)\in E(BS_n)$ if and only if there exists an integer $k$ with $k\in[2, n]$ such that $y_{k-1}=x_k$, $y_k=x_{k-1}$, and $y_i=x_i$ for every $i\in[1,n]-\{k-1,k\}$ or $y_1=x_k$, $y_k=x_1$, and $y_i=x_i$ for every $i\in[2,n]-\{k\}$.

\vskip.2cm

By Definition 2.1, $BS_n$ is a bipartite and $(2n-3)$-regular graph of order $n!$. Fig. 1 illustrates $BS_2$, $BS_3$, and $BS_4$,
respectively.

\begin{figure}[ht]
   \begin{center}
    \includegraphics[scale=0.8]{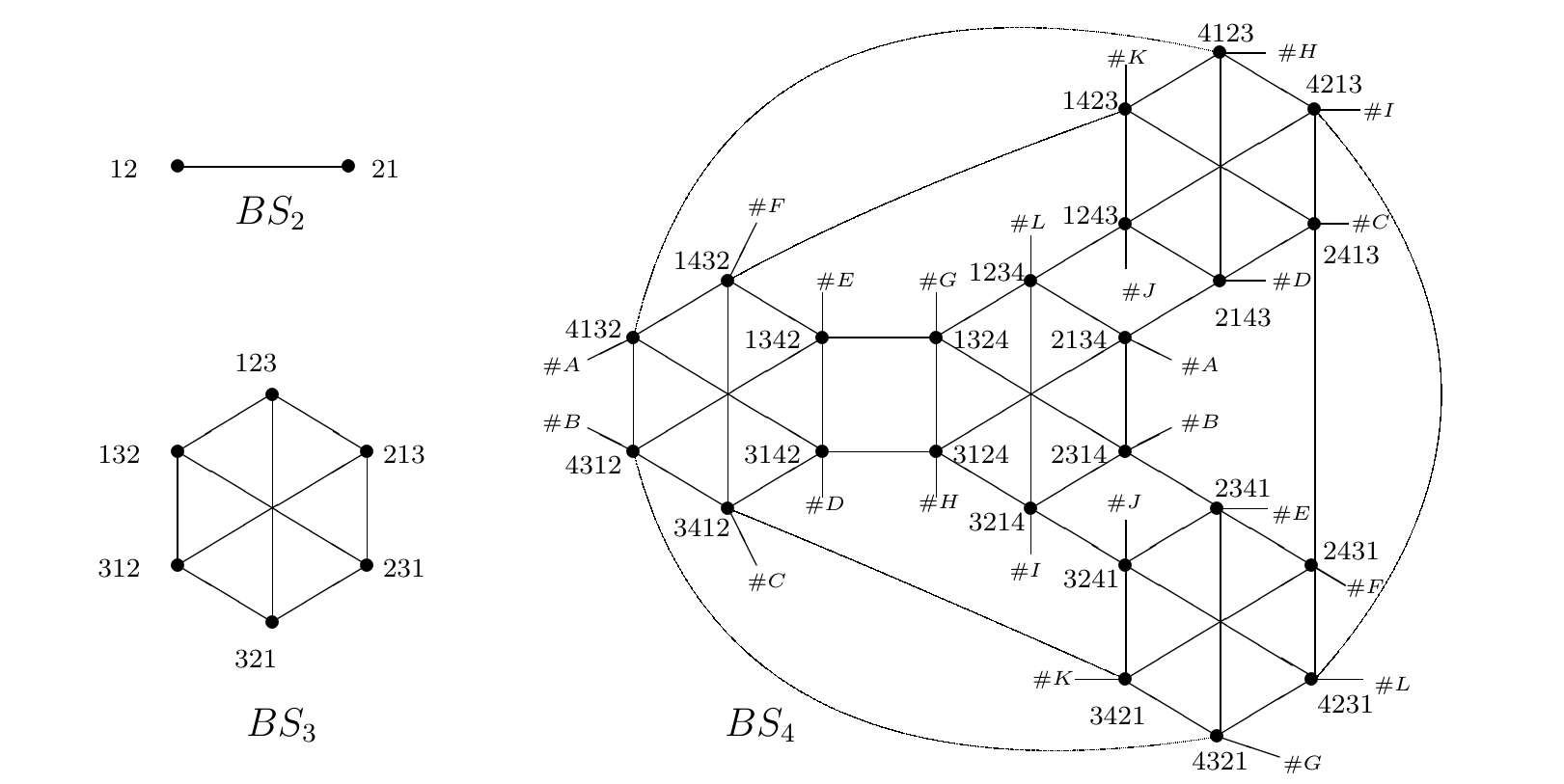}
    \end{center}
    \vspace{-0.5cm}\caption{\label{F2-5} Illustration of $BS_n$ for $n=2,3,4.$}
\end{figure}

\vskip.3cm

Let integers $j,k\in[1,n] $ with $j\neq k$. Let $x=x_1x_2\cdots x_n\in V(BS_n)$ and
``$\circ$" be an operation such that $y=y_1y_2\cdots y_n=x\circ(j,k)$ if and only if $x_{j}=y_{k}$, $x_{k}=y_{j}$, and $x_i=y_i$ for every $i\in[1,n]-\{j,k\}$. Thus $(x,y)\in E(BS_n)$ if and only if $y=x\circ(k-1,k)$ or $y=x\circ(1,k)$ for some $k\in[2, n]$. Let $x^-=x\circ(n-1,n)$ and $x^+=x\circ(1,n)$ for simplicity.
Let $BS_n^i$ be the induced subgraph of $BS_n$ by the vertex set
$V(BS_n^i)=\{x=x_1x_2\cdots x_n\in V(BS_n)~|~x_n=i\}$ for every $i\in[1, n]$. By Definition 2.1, $BS_n^i\cong BS_{n-1}$ for every $i\in[1,n]$. It is obvious that if $x\in V(BS_n^i)$, $x^-\in V(BS_n^j)$, and $x^+\in V(BS_n^k)$, then $i,j,k$ are three distinct integers in $[1,n]$.
Set $E_{i,j}(BS_n)=\{(x,y)\in E(BS_n)~|~x\in V(BS_n^i),~y\in V(BS_n^j)\}$ for any $i,j\in[1, n]$ with $i\neq j$. For any arbitrary edge set $F_e\subseteq E(BS_n)$, denote $F_e^i=F_e\cap E(BS_n^i)$ for every $i\in [1,n]$ and let $F_e^0=F_e-\cup_{i=1}^nF_e^i$.
For any $L\subseteq [1,n]$, let $BS_n^L$ be the subgraph of $BS_n$ induced by $\cup_{i\in L}V(BS_n^i)$.

Now we give some properties of $BS_n$.

\vskip.2cm

{\bf Lemma 2.2 \cite{Cai} } {\em Let $n$ be an integer with $n\geq3$. Then

 (1) $|E_{i,j}(BS_n)|=2(n-2)!$ for any $i,j\in[1,n]$ with $i\neq j$;

(2) $\{u^+,u^-\}\cap\{v^+,v^-\}=\emptyset$ for any $u,v\in V(BS_n^k)$ ($k\in[1,n]$) with $u\neq v$;

(3) $u^+\in V(BS_n^{[3,n]})$ or $u^-\in V(BS_n^{[3,n]})$ for any $u\in V(BS_n^{[1,2]})$.

%(4) For $u,v\in V(BS_n)$, $|N_{BS_n}(u)\cap N_{BS_n}(v)|\leq3$.
}

\vskip.2cm

{\bf Lemma 2.3 \cite{Wang3} } {\em $\lambda(BS_n)=2n-3$ for $n\geq3$.% and if $F_e$ is an edge cut of $BS_n$ with $|F_e|=2n-3$, then there exists a vertex $v$ such that $F_e=E_{BS_n}(v)$.
}

\vskip.2cm
{\bf Lemma 2.4 \cite{Wang3} } {\em Let $F_e\subseteq E(BS_n)$ with $|F_e|\leq4n-9$ for $n\geq3$. If $BS_n-F_e$ is disconnected, then $BS_n-F_e$ has two components, one of which is an isolated vertex.}

\vskip.2cm
{\bf Lemma 2.5 } {\em Let $F_e\subseteq E(BS_3)$ with $|F_e|\leq4$. If $BS_3-F_e$ is disconnected, then $BS_3-F_e$ has two components, one of which is an isolated vertex or an edge.}
\vskip.2cm

\p If $|F_e|\leq3$, then the lemma holds by Lemma 2.4. Now we consider the case that $|F_e|=4$ and $BS_3-F_e$ is disconnected. Let $H_1,H_2,\cdots,H_k$ be the $k$ components of $BS_3-F_e$ with $|V(H_1)|\geq|V(H_2)|\geq\cdots\geq|V(H_k)|$ and $k\geq2$. Since $|V(BS_3)|=3!=6$, $3\geq|V(H_2)|\geq\cdots\geq|V(H_k)|$. If $|V(H_2)|=3$, then $H_2=P_3$ as $BS_3$ is bipartite. Thus $|F_e|\geq 2\times2+1=5>4$, a contradiction. Hence $|V(H_2)|\leq2$. Now we claim that $k=2$.
Suppose, to the contrary, that $k\geq3$. Note that $BS_3$ is bipartite.
If $|V(H_2)|=|V(H_3)|=1$, then $|F_e|\geq 2\times3-1=5>4$, a contradiction. If $|V(H_2)|=|V(H_3)|=2$, then $|F_e|\geq 4\times2-2=6>4$, a contradiction. If $|V(H_2)|=2$ and $|V(H_3)|=1$, then $|F_e|\geq 2\times2+3-1=6>4$, a contradiction. Thus $k=2$ and the lemma holds.\q
\vskip.2cm

{\bf Lemma 2.6} {\em Let $F_e\subseteq E(BS_4)$ with $|F_e|\leq10$. If $BS_4-F_e$ is disconnected, then $BS_4-F_e$ has a component $H$ with $|V(H)|\geq 4!-2$.}
\vskip.2cm

\p Suppose that $BS_4-F_e$ is disconnected. Without loss of generality, we assume $|F_e^1|\geq|F_e^2|\geq|F_e^3|\geq|F_e^4|$. Since $n=4$, $|E_{i,j}(BS_4)|=2\times(4-2)!=4$ for $i,j\in[1,4]$ with $i\neq j$ by Lemma 2.2 (1).
Since $|F_e|\leq10$, $|F_e^4|\leq2$. Hence $BS_4^4-F_e^4$ is connected by Lemma 2.3. Let $H$ be the component of $BS_4-F_e$ containing $BS_4^4-F_e^4$ as a subgraph. Now we will consider the following three cases.

{\bf Case 1.} {\em $|F_e^1|\geq5$.}

In this case, $|F_e^4|\leq|F_e^3|\leq2$; otherwise $|F_e|\geq5+2\times3=11>10$, a contradiction. Thus $BS_4^3-F_e^3$ is connected by Lemma 2.3.

{\bf Subcase 1.1.} {\em $|F_e^2|\geq3$.}

In this subcase, $|F_e^0|\leq 10-5-3=2$. Since $|E_{3,4}(BS_4)-F_e|\geq |E_{3,4}(BS_4)|-|F_e^0|\geq4-2=2>0$, $BS_4^{[3,4]}-F_e$ is a subgraph of $H$. Since $|F_e^0|\leq2$, $|V(H)|\geq4!-2$ by Lemma 2.2 (3).

{\bf Subcase 1.2.} {\em $|F_e^2|\leq2$.}

In this subcase, $|F_e^0|\leq 10-5=5$ and $BS_4^i-F_e^i$ ($i=2,3,4$) is connected by Lemma 2.3. We claim that $E_{2,3}(BS_4)-F_e\neq\emptyset$ or $E_{2,4}(BS_4)-F_e\neq\emptyset$; otherwise $|F_e^0|\geq|E_{2,3}(BS_4)|+|E_{2,4}(BS_4)|=2\times4=8>5$, a contradiction. Without loss of generality, we assume $E_{2,3}(BS_4)-F_e\neq\emptyset$. Similarly, we can get
$E_{2,4}(BS_4)-F_e\neq\emptyset$ or $E_{3,4}(BS_4)-F_e\neq\emptyset$. Thus $BS_4^{[2,4]}-F_e$ is a subgraph of $H$. If $v\in V(BS_4^1)$, then $v^+\in V(BS_4^{[2,4]})$ and $v^-\in V(BS_4^{[2,4]})$. Since $|F_e^0|\leq5<2\times3$, $|V(H)|\geq4!-2$ by Lemma 2.2 (2).

{\bf Case 2.} {\em $3\leq|F_e^1|\leq4$.}

We will consider the following subcases.

{\bf Subcase 2.1.} {\em $|F_e^3|\geq3$.}

Since $3\leq|F_e^3|\leq|F_e^2|\leq|F_e^1|\leq4$ and $|F_e|\leq10$, we have $|F_e^3|=|F_e^2|=3$ and $|F_e^0|\leq10-3\times3=1$. Hence $BS_4^i-F_e^i$ has a component $H_i$ with $|V(H_i)|\geq3!-1$ for $i=2,3$ by Lemma 2.4. Since $|F_e^1|\leq4$, $BS_4^1-F_e^1$ has a component $H_1$ with $|V(H_1)|\geq3!-2$ by Lemma 2.5.
Since $|E_{BS_4}(V(H_i),V(BS_4^4))-F_e|\geq |E_{i,4}(BS_4)|-(3!-|V(H_i)|)-|F_e^0|\geq4-2-1>0$ for every $i\in[1,3]$, $H_i$ is a subgraph of $H$. If $BS_4^1-F_e^1$ is connected, then $|V(H)|\geq4!-2$. If $|V(H_1)|\geq3!-1$ and $BS_4^2-F_e^2$ or $BS_4^3-F_e^3$ is connected, then $|V(H)|\geq4!-2$. If $|V(H_1)|\geq3!-2$, both $BS_4^2-F_e^2$ and $BS_4^3-F_e^3$ are connected, then $|V(H)|\geq4!-2$.
Hence we just need to consider the following three conditions.

{\bf Subcase 2.1.1.} {\em $|V(H_1)|=|V(H_2)|=|V(H_3)|=3!-1$.}

Let $u_i\in V(BS_4^i)-V(H_i)$ for every $i\in[1,3]$. If $u_i\in V(H)$ for some $i\in[1,3]$, then the lemma holds. Now we suppose that $u_i\not\in V(H)$ for every $i\in[1,3]$. Note that $BS_4$ is bipartite.
If $u_1,u_2,u_3$ are three isolated vertices in $BS_4-F_e$, then $|F_e|\geq3\times5-2=13>10$, a contradiction. If $u_1,u_2,u_3$ form an edge and an isolated vertex in $BS_4-F_e$, then $|F_e|\geq2\times4+5-1=12>10$, a contradiction. If $u_1,u_2,u_3$ form a $P_3$ in $BS_4-F_e$, then $|F_e|\geq2\times4+3=11>10$, a contradiction.

{\bf Subcase 2.1.2.} {\em $|V(H_1)|=3!-2,~|V(H_2)|=|V(H_3)|=3!-1$.}

Let $u_i\in V(BS_4^i)-V(H_i)$ for $i=2,3$. Let $u_{11},u_{12}\in V(BS_4^1)-V(H_1)$ with $u_{11}\neq u_{12}$. Hence $|F_e^1|=4$, $|F_e^0|=0$, and $(u_{11}, u_{12})\in E(BS_4^1)-F_e$ by Lemmas 2.4 and 2.5. If $u_{11}\in V(H)$ or $u_{12}\in V(H)$, then the lemma holds. Now we suppose that $u_{11}\not\in V(H)$ and $u_{12}\not\in V(H)$. Hence $\{u_{11}^+,u_{11}^-\}=\{u_2,u_3\}$ as $|F_e^0|=0$. Thus $\{u_{12}^+,u_{12}^-\}\subseteq V(H)$ by Lemma 2.2 (2). Since $|F_e^0|=0$, $u_{12}\in V(H)$, a contradiction.

{\bf Subcase 2.1.3.} {\em $|V(H_1)|=3!-2,|V(H_2)|=3!-1,|V(H_3)|=3!$ or $|V(H_1)|=3!-2,|V(H_2)|=3!,|V(H_3)|=3!-1$.}

Without loss of generality, we assume $|V(H_1)|=3!-2,|V(H_2)|=3!-1,|V(H_3)|=3!$. Let $u_{11},u_{12}\in V(BS_4^1)-V(H_1)$ with $u_{11}\neq u_{12}$ and $u_2\in V(BS_4^2)-V(H_2)$. Hence $|F_e^1|=4$, $|F_e^0|=0$, and $(u_{11}, u_{12})\in E(BS_4^1)-F_e$ by Lemmas 2.4 and 2.5. Since $|F_e^0|=0$, $u_{11}^+\in V(H)$ or $u_{11}^-\in V(H)$. Hence $u_1\in V(H)$, the lemma holds.

{\bf Subcase 2.2.} {\em $|F_e^3|\leq2$.}

In this subcase, $|F_e^0|\leq10-3=7$. By Lemma 2.3, $BS_4^3-F_e^3$ is connected. Now we consider
the following three conditions.

{\bf Subcase 2.2.1.} {\em $|F_e^2|\leq2$.}

$BS_4^2-F_e^2$ is connected by Lemma 2.3. We claim that $E_{2,3}(BS_4)-F_e\neq\emptyset$ or $E_{2,4}(BS_4)-F_e\neq\emptyset$; otherwise $|F_e^0|\geq|E_{2,3}(BS_4)|+|E_{2,4}(BS_4)|=2\times4=8>7$, a contradiction. Without loss of generality, we assume $E_{2,3}(BS_4)-F_e\neq\emptyset$. Similarly, we can get  $E_{2,4}(BS_4)-F_e\neq\emptyset$ or $E_{3,4}(BS_4)-F_e\neq\emptyset$. Hence $BS_4^{[2,4]}-F_e$ is a subgraph of $H$. Since $3\leq|F_e^1|\leq4$, $BS_4^1-F_e^1$ has a component $H_1$ such that $|V(H_1)|\geq3!-2$ by Lemma 2.5. Since $\{u^+,u^-\}\subseteq V(BS_4^{[2,4]})$ for every $u\in V(BS_4^1)$, $|E_{BS_4}(V(H_1),V(BS_4^{[2,4]}))-F_e|\geq |E_{1,2}(BS_4)|+|E_{1,3}(BS_4)|+|E_{1,4}(BS_4)|-2|V(BS_4^1)-V(H_1)|-|F_e^0|\geq3\times4-2\times2-7>0$. Thus $H_1$ is a subgraph of $H$ and the lemma holds.

{\bf Subcase 2.2.2.} {\em $|F_e^2|=3$.}

In this subcase, we have $|F_e^0|\leq10-3-3=4$. If $BS_4^2-F_e^2$ is connected, then the lemma holds by the same argument as that of Subcase 2.2.1.
%We claim that $E_{2,3}(BS_4)-F_e\neq\emptyset$ or $E_{2,4}(BS_4)-F_e\neq\emptyset$. If $E_{2,3}(BS_4)-F_e=E_{2,4}(BS_4)-F_e=\emptyset$. Then $|F_e^0|\geq|E_{2,3}(BS_4)|+|E_{2,4}(BS_4)|=2\times2\times(4-2)!=8>4$, a contradiction. Without loss of generality, we assume $E_{2,3}(BS_4)-F_e\neq\emptyset$. Similarly, we can get  $E_{2,4}(BS_4)-F_e\neq\emptyset$ or $E_{3,4}(BS_4)-F_e\neq\emptyset$. Hence $BS_4^{[2,4]}-F_e^{[2,4]}$ is a subgraph of $H$. Since $3\leq|F_e^1|\leq4$, $BS_4^1-F_e^1$ has a component $H_1$ such that $|V(H_1)|\geq3!-2$ by Lemma 2.5. Since $\{u^+,u^-\}\subseteq V(BS_4^{[2,4]})$ for every $u\in V(BS_4^1)$, $|E_{BS_4}(V(H_1),V(BS_4^{[2,4]}-F_e))|\geq |E_{1,2}(BS_4)|+|E_{1,3}(BS_4)|+|E_{1,4}(BS_4)|-2|V(BS_4^1)-V(H)|-|F_e^0|=3\times2\times(4-2)!-2\times2-4=1>0$. Thus $H_1$ is a subgraph of $H$ and the lemma holds.

Now we suppose that $BS_4^2-F_e^2$ is disconnected. Then by Lemma 2.4, $BS_4^2-F_e^2$ has a component $H_2$ such that $|V(H_2)|=3!-1$. Let $u_2\in V(BS_4^2)-V(H_2)$. We claim that $E_{BS_4}(V(H_2),V(BS_4^3))-F_e\neq\emptyset$ or $E_{BS_4}(V(H_2),V(BS_4^4))-F_e\neq\emptyset$; otherwise $|F_e^0|\geq|E_{BS_4}(V(H_2),V(BS_4^3))|+|E_{BS_4}(V(H_2),V(BS_4^4))|\geq 4-1+4-1=6>4$, a contradiction. Without loss of generality, we assume $E_{BS_4}(V(H_2),V(BS_4^3))-F_e\neq\emptyset$. Similarly, we can get $E_{BS_4}(V(H_2),V(BS_4^4))-F_e\neq\emptyset$ or $E_{3,4}(BS_4)-F_e\neq\emptyset$. Hence both $H_2$ and $BS_4^{[3,4]}-F_e$ are subgraphs of $H$. Since $3\leq|F_e^1|\leq4$, $BS_4^1-F_e^1$ has a component $H_1$ such that $|V(H_1)|\geq3!-2$ by Lemma 2.5. If $|V(H_1)|\geq3!-1$, then $|E_{BS_4}(V(H_1),V(BS_4^{[3,4]}))-F_e|\geq |E_{1,3}(BS_4)|+|E_{1,4}(BS_4)|-2|V(BS_4^1)-V(H)|-|F_e^0|\geq2\times4-2\times1-4=2>0$, which implies $H_1$ is a subgraph of $H$ and the lemma holds. Now we consider that $|V(H_1)|=3!-2$. Hence $|F_e^1|=4$ by Lemmas 2.4 and 2.5. Thus $|F_e^0|\leq10-4-3=3$ and $|E_{BS_4}(V(H_1),V(BS_4^{[3,4]}))-F_e|\geq |E_{1,3}(BS_4)|+|E_{1,4}(BS_4)|-2|V(BS_4^1)-V(H)|-|F_e^0|\geq2\times4-2\times2-3=1>0$, which implies $H_1$ is a subgraph of $H$. Let $u_{11},u_{12}\in V(BS_4^1)-V(H_1)$ with $u_{11}\neq u_{12}$. Then the lemma holds by the same argument as that of Subcase 2.1.1.

{\bf Subcase 2.2.3.} {\em $|F_e^2|=4$.}

Since $|F_e^2|\leq|F_e^1|$, $|F_e^2|=|F_e^1|=4$ and $|F_e^0|\leq10-4-4=2$. Since $|E_{3,4}(BS_4)-F_e|\geq|E_{3,4}(BS_4)|-|F_e^0|\geq4-2=2>0$, $BS_4^{[3,4]}-F_e$ is a subgraph of $H$. Since $|F_e^0|\leq2$, the lemma holds by Lemma 2.2 (3).

{\bf Case 3.} {\em $|F_e^1|\leq2$.}

In this case, $BS_4^i-F_e^i$ ($i=1,2,3,4$) is connected by Lemma 2.3. Now we claim that $E_{1,k}(BS_4)-F_e\neq\emptyset$ for some $k\in[2,4]$; otherwise $|F_e|\geq |E_{1,2}(BS_4)|+|E_{1,3}(BS_4)|+|E_{1,4}(BS_4)|=3\times4=12>10$, a contradiction. Without loss of generality, we assume $E_{1,2}(BS_4)-F_e\neq\emptyset$. Suppose $E_{1,3}(BS_4)-F_e\neq\emptyset$ or $E_{2,3}(BS_4)-F_e\neq\emptyset$. Thus $BS_4^{[1,3]}-F_e$ is connected.
Similarly, we can get $E_{k,4}(BS_4)-F_e\neq\emptyset$ for some $k\in[1,3]$, which implies $H=BS_4-F_e$ is connected, a contradiction. Hence $E_{1,3}(BS_4)-F_e=\emptyset$ and $E_{2,3}(BS_4)-F_e=\emptyset$. Thus $|F_e\cap(E_{1,3}(BS_4)\cup E_{2,3}(BS_4))|=2\times4=8$. Hence $|E_{k,4}(BS_4)\cap F_e|\leq10-8=2$ and $|E_{k,4}(BS_4)- F_e|\geq4-2=2>0$ for every $k\in[1,3]$. Hence $H=BS_4-F_e$ is connected, a contradiction.\q

\vskip.2cm

{\bf Lemma 2.7} {\em Let $F_e\subseteq E(BS_n)$ with $|F_e|\leq6n-14$ for $n\geq3$. If $BS_n-F_e$ is disconnected, then $BS_n-F_e$ has a component $H$ with $|V(H)|\geq n!-2$.}
\vskip.2cm

\p We prove this lemma by induction on $n$. For $n=3,4$, the result holds by Lemmas 2.5 and 2.6. Assume $n\geq5$
and $BS_n-F_e$ is disconnected. Without loss of generality, we assume $|F_e^1|\geq|F_e^2|\geq\cdots\geq|F_e^n|$.
Since $|F_e|\leq6n-14$, $|F_e^n|\leq\cdots\leq|F_e^4|\leq 2n-6$; otherwise $|F_e|\geq4(2n-5)>6n-14$ for $n\geq5$, a contradiction. Hence $BS_n^i-F_e^i$ is connected for every $i\in[4,n]$ by Lemma 2.3. Let $H$ be the component of $BS_n-F_e$ containing $BS_n^n-F_e^n$ as a subgraph.
Now we will consider the following four cases.

{\bf Case 1.} {\em $|F_e^1|\geq6n-19$.}

In this case, $|F_e^0|\leq(6n-14)-(6n-19)=5$ and $|F_e^3|\leq2\leq2n-6$ for $n\geq5$. Hence $BS_n^3-F_e^3$ is connected by Lemma 2.3. Since $|E_{i,j}(BS_n)-F_e|\geq|E_{i,j}(BS_n)|-|F_e^0|\geq2(n-2)!-5>0$ for $i,j\in[3,n]$ with $i\neq j$ and $n\geq5$, $BS_n^{[3,n]}-F_e$ is a subgraph of $H$.

Suppose $BS_n^2-F_e^2$ is connected. Since $|E_{2,3}(BS_n)-F_e|\geq|E_{2,3}(BS_n)|-|F_e^0|\geq2(n-2)!-5>0$ for $n\geq5$, $BS_n^2-F_e^2$ is a subgraph of $H$. Note that $\{u^+,u^-\}\subseteq V(BS_n^{[2,n]})$ for every $u\in V(BS_n^1)$.
Since $|F_e^0|\leq5<2\times3$, we have $|V(H)|\geq n!-2$ by Lemma 2.2 (2).

Now we consider that $BS_n^2-F_e^2$ is disconnected. Then $2n-5\leq|F_e^2|\leq5$, which implies $n=5$, $|F_e^2|=5$, and $|F_e^0|=0$. Since $|F_e^0|=0$, $H=BS_n-F_e$ is connected by Lemma 2.2 (3), a contradiction.

{\bf Case 2.} {\em $4n-12\leq|F_e^1|\leq6n-20$.}

In this case, $|F_e^0|\leq(6n-14)-(4n-12)=2n-2$ and $|F_e^3|\leq2n-6$; otherwise $|F_e|\geq2(2n-5)+(4n-12)=8n-22>6n-14$ for $n\geq5$, a contradiction. Thus $BS_n^i-F_e^i$ is connected for every $i\in[3,n]$ by Lemma 2.3. Since $|E_{i,j}(BS_n)-F_e|\geq|E_{i,j}(BS_n)|-|F_e^0|\geq2(n-2)!-(2n-2)>0$ for $i,j\in[3,n]$ with $i\neq j$ and $n\geq5$, $BS_n^{[3,n]}-F_e$ is a subgraph of $H$.

Suppose $BS_n^2-F_e^2$ is connected. Since $|E_{2,3}(BS_n)-F_e|\geq|E_{2,3}(BS_n)|-|F_e^0|\geq2(n-2)!-(2n-2)>0$ for $n\geq5$, $BS_n^2-F_e^2$ is a subgraph of $H$. Since $4n-12\leq|F_e^1|\leq6n-20$, $BS_n^1-F_e^1$ has a component $H_1$ with $|V(H_1)|\geq(n-1)!-2$ by induction hypothesis. Since $|E_{BS_n}(V(H_1),V(BS_n^2))-F_e|\geq|E_{1,2}(BS_n)|-|V(BS_n^1)-V(H_1)|-|F_e^0|\geq 2(n-2)!-2-(2n-2)>0$ for $n\geq5$, $H_1$ is a subgraph of $H$. Thus $|V(H)|\geq n!-2$.

Now we consider that $BS_n^2-F_e^2$ is disconnected. Hence $2n-5\leq|F_e^2|\leq|F_e^1|\leq6n-20$ and $|F_e^0|\leq(6n-14)-(4n-12)-(2n-5)=3$.
%Hence  $BS_n^i-F_e^i$ has a component $H_i$ with $|V(H_i)|\geq(n-1)!-2$ for $i=1,2$ by induction hypothesis. Since $|E_{BS_n}(V(H_i),V(BS_n^3))-F_e|\geq|E_{i,3}(BS_n)|-|V(BS_n^i)-V(H_i)|-|F_e^0|\geq 2(n-2)!-2-3>0$ for $i=1,2$ and $n\geq5$, $H_i$ is a subgraph of $H$.
Since $|F_e^0|\leq3$, $|V(BS_n)-V(H)|\leq3$ by Lemma 2.2 (3). If $|V(BS_n)-V(H)|\leq2$, then the lemma holds. Now we suppose $|V(BS_n)-V(H)|=3$ and $V(BS_n)-V(H)=\{u_1,u_2,u_3\}$. Note that $BS_n$ is bipartite. If $u_1,u_2,u_3$ are three isolated vertices in $BS_n-F_e$, then $|F_e|\geq3(2n-3)-2=6n-11>6n-14$, a contradiction. If $u_1,u_2,u_3$ form an edge and an isolated vertex in $BS_n-F_e$, then $|F_e|\geq2(2n-4)+(2n-3)-1=6n-12>6n-14$, a contradiction. If $u_1,u_2,u_3$ form a $P_3$ in $BS_n-F_e$, then $|F_e|\geq2(2n-4)+(2n-5)=6n-13>6n-14$, a contradiction.

{\bf Case 3.} {\em $2n-5\leq|F_e^1|\leq4n-13$.}

In this case, $|F_e^0|\leq(6n-14)-(2n-5)=4n-9$.

{\bf Subcase 3.1.} {\em $|F_e^2|\leq2n-6$.}

In this subcase, $BS_n^i-F_e^i$ is connected for every $i\in[2,n]$ by Lemma 2.3. Since $|E_{i,j}(BS_n)-F_e|\geq|E_{i,j}(BS_n)|-|F_e^0|\geq2(n-2)!-(4n-9)>0$ for $i,j\in[2,n]$ with $i\neq j$ and $n\geq5$, $BS_n^{[2,n]}-F_e$ is a subgraph of $H$. Since $2n-5\leq|F_e^1|\leq4n-13$, $BS_n^1-F_e^1$ has a component $H_1$ with $|V(H_1)|\geq(n-1)!-1$ by Lemma 2.4. Since $|E_{BS_n}(V(H_1),V(BS_n^{[2,3]}))-F_e|\geq|E_{1,2}(BS_n)|+|E_{1,3}(BS_n)|-2|V(BS_n^1)-V(H_1)|-|F_e^0|\geq2\times2(n-2)!-2\times1-(4n-9)>0$ for $n\geq5$, $H_1$ is a subgraph of $H$ and $|V(H)|\geq n!-1$.

{\bf Subcase 3.2.} {\em $2n-5\leq|F_e^2|\leq4n-13$.}

In this subcase, $|F_e^0|\leq(6n-14)-2(2n-5)=2n-4$. If $|F_e^3|\leq2n-6$, then $BS_n^i-F_e^i$ is connected for every $i\in[3,n]$ by Lemma 2.3. Since $|E_{i,j}(BS_n)-F_e|\geq|E_{i,j}(BS_n)|-|F_e^0|\geq2(n-2)!-(2n-4)>0$ for $i,j\in[3,n]$ with $i\neq j$ and $n\geq5$, $BS_n^{[3,n]}-F_e$ is a subgraph of $H$. Since $2n-5\leq|F_e^2|\leq|F_e^1|\leq4n-13$, $BS_n^k-F_e^k$ has a component $H_k$ with $|V(H_k)|\geq(n-1)!-1$ for $k=1,2$ by Lemma 2.4. Since $|E_{BS_n}(V(H_k),V(BS_n^3))-F_e|\geq|E_{k,3}(BS_n)|-|V(BS_n^k)-V(H_k)|-|F_e^0|\geq2(n-2)!-1-(2n-4)>0$ for $k\in[1,2]$ and $n\geq5$, both $H_1$ and $H_2$ are subgraphs of $H$. Thus $|V(H)|\geq n!-2$.

Suppose $|F_e^3|\geq2n-5$. Then $|F_e^0|\leq(6n-14)-3(2n-5)=1$. Since $|E_{i,j}(BS_n)-F_e|\geq|E_{i,j}(BS_n)|-|F_e^0|\geq2(n-2)!-1>0$ for $i,j\in[4,n]$ with $i\neq j$ and $n\geq5$, $BS_n^{[4,n]}-F_e$ is a subgraph of $H$. Since $2n-5\leq|F_e^3|\leq|F_e^2|\leq|F_e^1|\leq4n-13$, $BS_n^k-F_e^k$ has a component $H_k$ with $|V(H_k)|\geq(n-1)!-1$ for every $k\in[1,3]$ by Lemma 2.4. Since $|E_{BS_n}(V(H_k),V(BS_n^4))-F_e|\geq|E_{k,4}(BS_n)|-|V(BS_n^k)-V(H_k)|-|F_e^0|\geq2(n-2)!-1-1>0$ for $k\in[1,3]$ and $n\geq5$, $H_i$ is a subgraph of $H$ for every $k\in[1,3]$. If $BS_n^k-F_e^k$ is connected for some $k\in[1,3]$, then $|V(H)|\geq n!-2$. Now we consider that $|V(H_1)|=|V(H_2)|=|V(H_3)|=(n-1)!-1$. Let $u_k\in V(BS_n^k)-V(H_k)$ for every $k\in[1,3]$. Then the lemma holds by the same argument as that of Case 2.
%Note that $BS_n$ is a bipartite graph. If $u_1,u_2,u_3$ are three isolated vertices in $BS_n-F_e$, then $|F_e|\geq3(2n-3)-2=6n-11>6n-14$, a contradiction. If $u_1,u_2,u_3$ form an edge and a isolated vertex in $BS_n-F_e$, then $|F_e|\geq2(2n-4)+(2n-3)-1=6n-12>6n-14$, a contradiction. If $u_1,u_2,u_3$ form a $P_3$ in $BS_n-F_e$, then $|F_e|\geq2(2n-4)+(2n-5)=6n-13>6n-14$, a contradiction.

{\bf Case 4.} {\em $|F_e^1|\leq2n-6$.}

In this case, $BS_n^i-F_e^i$ is connected for every $i\in[1,n]$ by Lemma 2.3. We claim that $E_{1,2}(BS_n)-F_e\neq\emptyset$ or $E_{1,3}(BS_n)-F_e\neq\emptyset$; otherwise $|F_e|\geq |E_{1,2}(BS_n)|+|E_{1,3}(BS_n)|=2\times2(n-2)!>6n-14$ for $n\geq5$, a contradiction. Without loss of generality, we assume $E_{1,2}(BS_n)-F_e\neq\emptyset$. Similarly, we can get $E_{1,i}(BS_n)-F_e\neq\emptyset$ or $E_{2,i}(BS_n)-F_e\neq\emptyset$ for every $i\in[3,n]$. Thus $H=BS_n-F_e$ is connected, a contradiction.\q

\vskip.2cm

{\bf Lemma 2.8} {\em Let $F_e\subseteq E(BS_4)$ with $|F_e|\leq11$. If $BS_4-F_e$ is disconnected, then $BS_4-F_e$ has a component $H$ with $|V(H)|\geq 4!-3$.}
\vskip.2cm

\p Suppose that $BS_4-F_e$ is disconnected. Without loss of generality, we assume $|F_e^1|\geq|F_e^2|\geq|F_e^3|\geq|F_e^4|$. Since $n=4$, $|E_{i,j}(BS_4)|=2\times(4-2)!=4$ for $i,j\in[1,4]$ with $i\neq j$ by Lemma 2.2 (1).
Since $|F_e|\leq11$, $|F_e^4|\leq2$. Hence $BS_4^4-F_e^4$ is connected by Lemma 2.3. Let $H$ be the component of $BS_4-F_e$ containing $BS_4^4-F_e^4$ as a subgraph. If $|F_e^1|\leq2$, then the lemma holds by the same argument as that of Case 3 of Lemma 2.6.
Hence we just consider the following two cases.

{\bf Case 1.} {\em $|F_e^1|\geq5$.}

Suppose that $|F_e^3|\geq3$. Since $|F_e^3|\leq|F_e^2|\leq|F_e^1|$, we have $|F_e^3|=|F_e^2|=3$, $|F_e^1|=5$, and $|F_e^0|=0$. Hence $BS_4^i-F_e^i$ has a component $H_i$ with $|V(H_i)|\geq3!-1$ for $i=2,3$ by Lemma 2.4. Since $|E_{BS_4}(V(H_i),V(BS_4^4))-F_e|\geq |E_{i,4}(BS_4)|-|V(BS_4^i)-V(H_i)|-|F_e^0|\geq4-1=3>0$ for $i=2,3$, both $H_2$ and $H_3$ are subgraphs of $H$. If $BS_4^3-F_e^3$ is a subgraph of $H$, then $H=BS_4-F_e$ is connected by Lemma 2.2 (3), a contradiction. Thus $|V(H_3)|=3!-1$ and there exists a vertex $u_3\in V(BS_4^3)-V(H)$. Since $|F_e^0|=0$ and $u_3\not\in V(H)$, $\{u_3^+,u_3^-\}\subseteq V(BS_4^{[1,2]})-V(H)$ and $|V(H_2)|=3!-1$. Let $\{u_3^+,u_3^-\}\cap V(BS_4^i)=u_i$ for $i=1,2$. Since $BS_4$ is bipartite and $|V(H_2)|=|V(H_3)|=3!-1$, $\{u_1^+,u_1^-\}\cap V(H)\neq\emptyset$. Since $|F_e^0|=0$, $u_1\in V(H)$, which implies $u_3\in V(H)$, a contradiction.

Now we suppose that $|F_e^3|\leq2$. Then $BS_4^3-F_e^3$ is connected by Lemma 2.3. Hence $|V(H)|\geq 4!-3$ by the same argument as that of Case 1 of Lemma 2.6

%{\bf Subcase 1.1.} {\em $|F_e^2|\geq3$.}

%In this subcase, $|F_e^0|\leq 11-5-3=3$. Since $|E_{3,4}(BS_4)-F_e|\geq |E_{3,4}(BS_4)|-|F_e^0|\geq2\times(4-2)!-3=1>0$ by Lemma 2.2 (1), $BS_4^{[3,4]}-F_e^{[3,4]}$ is a subgraph of $H$. Since $|F_e^0|\leq3$, $|V(H)|\geq4!-3$ by Lemma 2.2 (3).

%{\bf Subcase 1.2.} {\em $0\leq|F_e^2|\leq2$.}

%In this subcase, $|F_e^0|\leq 11-5=6$ and $BS_4^i-F_e^i$ ($i=2,3,4$) is connected by Lemma 2.3. We claim that $E_{2,3}(BS_4)-F_e\neq\emptyset$ or $E_{2,4}(BS_4)-F_e\neq\emptyset$; otherwise $|F_e^0|\geq|E_{2,3}(BS_4)|+|E_{2,4}(BS_4)|=2\times2\times(4-2)!=8>6$, a contradiction. Without loss of generality, we assume $E_{2,3}(BS_4)-F_e\neq\emptyset$. Similarly, we can get $E_{2,4}(BS_4)-F_e\neq\emptyset$ or $E_{3,4}(BS_4)-F_e\neq\emptyset$. Thus $BS_4^{[2,4]}-F_e^{[2,4]}$ is a subgraph of $H$. If $v\in V(BS_4^1)$, then $v^+\in V(BS_4^{[2,4]})$ and $v^-\in V(BS_4^{[2,4]})$. Since $|F_e^0|\leq6<2\times4$, $|V(H)|\geq4!-3$ by Lemma 2.2 (2).

{\bf Case 2.} {\em $3\leq|F_e^1|\leq4$.}

We will consider the following subcases.

{\bf Subcase 2.1.} {\em $|F_e^3|\geq3$.}

Since $3\leq|F_e^3|\leq|F_e^2|\leq|F_e^1|\leq4$ and $|F_e|\leq11$, we have $|F_e^3|=3$. Hence $BS_4^3-F_e^3$ has a component $H_3$ such that $|V(H_3)|\geq3!-1$ by Lemma 2.4.

{\bf Subcase 2.1.1.} {\em $|F_e^2|=4$.}

In this subcase, $|F_e^1|=4$ and $|F_e^0|=0$.
By Lemma 2.5, $BS_4^i-F_e^i$ has a component $H_i$ such that $|V(H_i)|\geq3!-2$ for $i=1,2$.
Since $|E_{BS_4}(V(H_i),V(BS_4^4))-F_e|\geq |E_{i,4}(BS_4)|-(3!-|V(H_i)|)-|F_e^0|\geq4-2-0>0$ for $i\in[1,3]$, $H_i$ is a subgraph of $H$ for every $i\in[1,3]$. If $BS_4^3-F_e^3$ is a subgraph of $H$, then $H=BS_4-F_e$ by Lemma 2.2 (3), a contradiction. Hence $|V(H_3)|=3!-1$ and there exists a vertex $u_3\in V(BS_4^3)-V(H)$. Since $|F_e^0|=0$ and $u_3\not\in V(H)$, $\{u_3^+,u_3^-\}\subseteq V(BS_4^{[1,2]})-V(H)$. Let $\{u_3^+,u_3^-\}\cap V(BS_4^i)=u_i$ for $i=1,2$. Since $BS_4$ is bipartite and $|V(H_3)|=3!-1$, there exists a vertex $u_2'\in V(BS_4^2)-V(H)-\{u_2\}$ such that $(u_1,u_2')\in E(BS_4)$. Thus $|V(H_2)|=3!-2$ and $(u_2,u_2')\in E(BS_4^2)-F_e$ by Lemma 2.5. Similarly, there exists a vertex $u_1'\in V(BS_4^1)-V(H)-\{u_1\}$ such that $(u_1',u_2)\in E(BS_4)$, $|V(H_1)|=3!-2$, and $(u_1,u_1')\in E(BS_4^1)-F_e$. Since $|V(H_3)|=3!-1$ and $BS_4$ is bipartite, $\{u_1'^+,u_1'^-\}-\{u_2\}\subseteq V(H)$ by Lemma 2.2 (3). Since $|F_e^0|=0$, $u_1'\in V(H)$, which implies $u_2\in V(H)$, a contradiction.

{\bf Subcase 2.1.2.} {\em $|F_e^2|=3$.}

By Lemma 2.4, $BS_4^2-F_e^2$ has a component $H_2$ such that $|V(H_2)|\geq3!-1$.

Suppose $|F_e^1|=3$, then $|F_e^0|\leq11-3\times3=2$. By Lemma 2.4, $BS_4^1-F_e^1$ has a component $H_1$ such that $|V(H_1)|\geq3!-1$.
Since $|E_{BS_4}(V(H_i),V(BS_4^4))-F_e|\geq |E_{i,4}(BS_4)|-(3!-|V(H_i)|)-|F_e^0|\geq4-1-2>0$ for $i\in[1,3]$, $H_i$ is a subgraph of $H$ for every $i\in[1,3]$. Thus $|V(H)|\geq 4!-3$.

Suppose $|F_e^1|=4$, then $|F_e^0|\leq11-4-2\times3=1$. By Lemma 2.5, $BS_4^1-F_e^1$ has a component $H_1$ such that $|V(H_1)|\geq3!-2$. Since $|E_{BS_4}(V(H_i),V(BS_4^4))-F_e|\geq |E_{i,4}(BS_4)|-(3!-|V(H_i)|)-|F_e^0|\geq4-2-1>0$ for $i\in[1,3]$, $H_i$ is a subgraph of $H$ for every $i\in[1,3]$. If $|V(H_1)|\geq3!-1$, then $|V(H)|\geq4!-3$. If $|V(H_2)|=3!$ or $|V(H_3)|=3!$, then $|V(H)|\geq4!-3$. Now we consider that $|V(H_1)|=3!-2$ and $|V(H_2)|=|V(H_3)|=3!-1$. Let $\{u_{11},u_{12}\}\subseteq V(BS_4^1)-V(H_1)$ with $u_{11}\neq u_{12}$. Then $(u_{11}, u_{12})\in E(BS_4^1)-F_e$ by Lemma 2.5.
If $u_{11}\in V(H)$ or $u_{12}\in V(H)$, then $|V(H)|\geq4!-2$. We suppose that $u_{11}\not\in V(H)$ and $u_{12}\not\in V(H)$. Since $BS_4$ is bipartite, $|V(H_2)|=|V(H_3)|=3!-1$, and $|F_e^0|\leq1$, there exists a vertex  $v\in\{u_{11}^+,u_{11}^-,u_{12}^+,u_{12}^-\}\cap V(H)$ such that $(u_{11},v)\in E(BS_4)-F_e$ or $(u_{12},v)\in E(BS_4)-F_e$ by Lemma 2.2 (2), which implies $u_{11}\in V(H)$ and $u_{12}\in V(H)$, a contradiction.

{\bf Subcase 2.2.} {\em $ |F_e^3|\leq2$.}

In this subcase, $|F_e^0|\leq11-3=8$. By Lemma 2.3, $BS_4^3-F_e^3$ is connected. If $|F_e^2|=4$, then the lemma holds by the same argument as that of Subcase 2.2.3 of Lemma 2.6.
Hence we just consider the following two conditions.

{\bf Subcase 2.2.1.} {\em $|F_e^2|\leq2$.}

By Lemma 2.3, $BS_4^2-F_e^2$ is connected.

Suppose $BS_4^{[2,4]}-F_e$ is connected. By Lemma 2.5, $BS_4^1-F_e^1$ has a component $H_1$ such that $|V(H_1)|\geq3!-2$. If $|V(H_1)|\geq3!-1$, then $|E_{BS_4}(V(H_1),V(BS_4^{[2,4]}))-F_e|\geq |E_{1,2}(BS_4)|+|E_{1,3}(BS_4)|+|E_{1,4}(BS_4)|-2(3!-|V(H_1)|)-|F_e^0|\geq3\times4-2\times1-8>0$. Hence $H_1$ is a subgraph of $H$ and $|V(H)|\geq4!-1$. Now we consider that $|V(H_1)|=3!-2$, which implies $|F_e^1|=4$ by Lemmas 2.4 and 2.5. Thus $|F_e^0|\leq11-4=7$ and $|E_{BS_4}(V(H_1),V(BS_4^{[2,4]}))-F_e|\geq |E_{1,2}(BS_4)|+|E_{1,3}(BS_4)|+|E_{1,4}(BS_4)|-2(3!-|V(H_1)|)-|F_e^0|\geq3\times4-2\times2-7>0$. Hence $H_1$ is a subgraph of $H$ and $|V(H)|\geq4!-2$.

Now we suppose that $BS_4^{[2,4]}-F_e$ is disconnected. Without loss of generality, we assume $E_{2,3}(BS_4)-F_e=E_{2,4}(BS_4)-F_e=\emptyset$. Hence $|F_e^0|\geq|E_{2,3}(BS_4)|+|E_{2,4}(BS_4)|=2\times4=8$. Since $|F_e|\leq11$ and $3\leq|F_e^1|\leq4$, we have $|F_e^1|=3$, $|F_e^2|=0$, and $F_e^0=E_{2,3}(BS_4)\cup E_{2,4}(BS_4)$. Thus $E_{3,4}(BS_4)-F_e=E_{3,4}(BS_4)$ and $BS_4^{[3,4]}-F_e$ is connected. By Lemma 2.4, $BS_4^1-F_e^1$ has a component $H_1$ such that $|V(H_1)|\geq3!-1$. Since $|E_{BS_4}(V(H_1),V(BS_4^3))-F_e|\geq |E_{1,3}(BS_4)|-(3!-|V(H_1)|)\geq4-1>0$, $H_1$ is a subgraph of $H$. Since $|E_{BS_4}(V(H_1),V(BS_4^2))-F_e|\geq |E_{1,2}(BS_4)|-(3!-|V(H_1)|)\geq4-1>0$, $BS_4^2-F_e^2$ is a subgraph of $H$. Thus $|V(H)|\geq4!-1$.

{\bf Subcase 2.2.2.} {\em $|F_e^2|=3$.}

In this subcase, we have $|F_e^0|\leq11-3-3=5$. By Lemma 2.4, $BS_4^2-F_e^2$ has a component $H_2$ such that $|V(H_2)|\geq3!-1$. By Lemma 2.5, $BS_4^1-F_e^1$ has a component $H_1$ such that $|V(H_1)|\geq3!-2$.

Suppose $BS_4^{[3,4]}-F_e$ is connected. Since $|E_{BS_4}(V(H_2),V(BS_4^{[3,4]}))-F_e|\geq |E_{2,3}(BS_4)|+|E_{2,4}(BS_4)|-2(3!-|V(H_2)|)-|F_e^0|\geq2\times4-2\times1-5>0$, $H_2$ is a subgraph of $H$.
Since $|E_{BS_4}(V(H_1),V(BS_4^{[3,4]})\cup V(H_2))-F_e|\geq |E_{1,3}(BS_4)|+|E_{1,4}(BS_4)|+|E_{1,2}(BS_4)|-2(3!-|V(H_1)|)-(3!-|V(H_2)|)-|F_e^0|\geq3\times4-2\times2-1-5>0$, $H_1$ is a subgraph of $H$ and $|V(H)|\geq4!-3$.

Now we suppose that $BS_4^{[3,4]}-F_e$ is disconnected. Then $|F_e\cap E_{3,4}(BS_4)|=|E_{3,4}(BS_4)|=4$ and $|F_e^0- E_{3,4}(BS_4)|\leq 11-3-3-4=1$. Since $|E_{BS_4}(V(H_2),V(BS_4^i))-F_e|\geq |E_{2,i}(BS_4)|-(3!-|V(H_2)|)-|F_e^0-E_{3,4}(BS_4)|\geq 4-1-1>0$ for $i=3,4$, both $H_2$ and $BS_4^i-F_e^i$ are subgraphs of $H$. Since $|E_{BS_4}(V(H_1),V(BS_4^3))-F_e|\geq |E_{1,3}(BS_4)|-(3!-|V(H_1)|)-|F_e^0-E_{3,4}(BS_4)|\geq 4-2-1>0$, $H_1$ is a subgraph of $H$. Thus $|V(H)|\geq4!-3$.\q

%{\bf Subcase 2.2.3.} {\em $|F_e^2|=4$.}

%Since $|F_e^2|\leq|F_e^1|$, $|F_e^2|=|F_e^1|=4$ and $|F_e^0|\leq11-4-4=3$. Since $|E_{3,4}(BS_4)-F_e|\geq|E_{3,4}(BS_4)|-|F_e^0|\geq2\times(4-2)!-3=1>0$, $BS_4^{[3,4]}-F_e^{[3,4]}$ is a subgraph of $H$. Since $|F_e^0|\leq3$, the lemma holds by Lemma 2.2 (3).

%{\bf Case 3.} {\em $0\leq|F_e^1|\leq2$.}

%The lemma holds by the same argument as that of Case 3 of Lemma 2.6.\q

%In this case, $BS_4^i-F_e^i$ ($i=1,2,3,4$) is connected by Lemma 2.3. Now we claim that $E_{1,k}(BS_4)-F_e\neq\emptyset$ for some $k\in[2,4]$; otherwise $|F_e|\geq |E_{1,2}(BS_4)|+|E_{1,3}(BS_4)|+|E_{1,4}(BS_4)|=3\times2\times(4-2)!=12>11$, a contradiction. Without loss of generality, we assume $E_{1,2}(BS_4)-F_e\neq\emptyset$. Suppose $E_{1,3}(BS_4)-F_e\neq\emptyset$ or $E_{2,3}(BS_4)-F_e\neq\emptyset$. Thus $BS_4^{[1,3]}-F_e^{[1,3]}$ is connected. Similarly, we can get $E_{k,4}(BS_4)-F_e\neq\emptyset$ for some $k\in[1,3]$, which implies $BS_4-F_e$ is connected, a contradiction. Hence $E_{1,3}(BS_4)-F_e=\emptyset$ and $E_{2,3}(BS_4)-F_e=\emptyset$. Thus $|F_e\cap(E_{1,3}(BS_4)\cup E_{2,3}(BS_4))|=2\times2\times(4-2)!=8$. Hence $|E_{k,4}(BS_4)\cap F_e|\leq11-8=3$, $|E_{k,4}(BS_4)- F_e|\geq2\times(4-2)!-3=1>0$ for every $k\in[1,3]$. Hence $BS_4-F_e$ is connected, a contradiction.
\vskip.2cm

{\bf Lemma 2.9} {\em Let $F_e\subseteq E(BS_n)$ with $|F_e|\leq8n-21$ for $n\geq3$. If $BS_n-F_e$ is disconnected, then $BS_n-F_e$ has a component $H$ with $|V(H)|\geq n!-3$.}
\vskip.2cm

\p We prove this lemma by induction on $n$. For $n=3,4$, the result holds by Lemmas 2.4 and 2.8. Assume $n\geq5$
and $BS_n-F_e$ is disconnected. Without loss of generality, we assume $|F_e^1|\geq|F_e^2|\geq\cdots\geq|F_e^n|$.
Since $|F_e|\leq8n-21$, $|F_e^n|\leq\cdots\leq|F_e^4|\leq 2n-6$; otherwise $|F_e|\geq4(2n-5)>8n-21$ for $n\geq5$, a contradiction. Hence $BS_n^i-F_e^i$ is connected for every $i\in[4,n]$ by Lemma 2.3. Let $H$ be the component of $BS_n-F_e$ containing $BS_n^n-F_e^n$ as a subgraph. If $|F_e^1|\leq2n-6$, then the lemma holds by the same argument as that of Case 4 of Lemma 2.7.
Now we will consider the following four cases.

{\bf Case 1.} {\em $|F_e^1|\geq8n-28$.}

In this case, $|F_e^0|\leq(8n-21)-(8n-28)=7$ and $|F_e^3|\leq4\leq2n-6$ for $n\geq5$. Thus $BS_n^i-F_e^i$ is connected for every $i\in[3,n]$ by Lemma 2.3. Since $|E_{i,j}(BS_n)-F_e|\geq|E_{i,j}(BS_n)|-|F_e^0|\geq2(n-2)!-7>0$ for $i,j\in[3,n]$ with $i\neq j$ and $n\geq5$, $BS_n^{[3,n]}-F_e$ is a subgraph of $H$.

Suppose $BS_n^2-F_e^2$ is connected. Since $|E_{2,3}(BS_n)-F_e|\geq|E_{2,3}(BS_n)|-|F_e^0|\geq2(n-2)!-7>0$ for $n\geq5$, $BS_n^2-F_e^2$ is a subgraph of $H$. Note that $\{u^+,u^-\}\subseteq V(BS_n^{[2,n]})$ for every $u\in V(BS_n^1)$.
Since $|F_e^0|\leq7<2\times4$, we have $|V(H)|\geq n!-3$ by Lemma 2.2 (2).

Now we consider that $BS_n^2-F_e^2$ is disconnected. Then $2n-5\leq|F_e^2|\leq7$ for $n\geq5$, which implies $5\leq|F_e^2|\leq4n-13$ and $|F_e^0|\leq(8n-21)-(8n-28)-5=2$. Since $|F_e^0|\leq2$, $|V(H)|\geq n!-2$ by Lemma 2.2 (3).

{\bf Case 2.} {\em $6n-19\leq|F_e^1|\leq8n-29$.}

In this case, $|F_e^0|\leq(8n-21)-(6n-19)=2n-2$ and $|F_e^3|\leq2n-6$; otherwise $|F_e|\geq2(2n-5)+(6n-19)=10n-29>8n-21$ for $n\geq5$, a contradiction. Thus $BS_n^i-F_e^i$ is connected for every $i\in[3,n]$ by Lemma 2.3. Since $|E_{i,j}(BS_n)-F_e|\geq|E_{i,j}(BS_n)|-|F_e^0|\geq2(n-2)!-(2n-2)>0$ for $i,j\in[3,n]$ with $i\neq j$ and $n\geq5$, $BS_n^{[3,n]}-F_e$ is a subgraph of $H$.

Suppose $BS_n^2-F_e^2$ is connected. Since $|E_{2,3}(BS_n)-F_e|\geq|E_{2,3}(BS_n)|-|F_e^0|\geq2(n-2)!-(2n-2)>0$ for $n\geq5$, $BS_n^2-F_e^2$ is a subgraph of $H$. Since $|F_e^1|\leq8n-29$, $BS_n^1-F_e^1$ has a component $H_1$ with $|V(H_1)|\geq(n-1)!-3$ by induction hypothesis. Since $|E_{BS_n}(V(H_1),V(BS_n^2))-F_e|\geq|E_{1,2}(BS_n)|-|V(BS_n^1)-V(H_1)|-|F_e^0|\geq 2(n-2)!-3-(2n-2)>0$ for $n\geq5$, $H_1$ is a subgraph of $H$. Hence $|V(H)|\geq n!-3$.

Now we suppose $BS_n^2-F_e^2$ is disconnected. Hence $2n-5\leq|F_e^2|\leq2n-2$ and $|F_e^0|\leq(8n-21)-(6n-19)-(2n-5)=3$. Since $|F_e^0|\leq3$, $|V(H)|\geq n!-3$ by Lemma 2.2 (3).

{\bf Case 3.} {\em $4n-12\leq|F_e^1|\leq6n-20$.}

In this case, $|F_e^0|\leq(8n-21)-(4n-12)=4n-9$. Since $|E_{i,j}(BS_n)-F_e|\geq|E_{i,j}(BS_n)|-|F_e^0|\geq2(n-2)!-(4n-9)>0$ for $i,j\in[4,n]$ with $i\neq j$ and $n\geq5$, $BS_n^{[4,n]}-F_e$ is a subgraph of $H$. Since $4n-12\leq|F_e^1|\leq6n-20$, $BS_n^1-F_e^1$ has a component $H_1$ with $|V(H_1)|\geq(n-1)!-2$ by Lemma 2.7.

{\bf Subcase 3.1.} {\em $4n-12\leq|F_e^2|\leq6n-20$.}

In this subcase, $|F_e^0|\leq(8n-21)-2(4n-12)=3$ and $|F_e^3|\leq3\leq2n-6$ for $n\geq5$. Hence $BS_n^3-F_e^3$ is connected by Lemma 2.3.
Since $|E_{3,4}(BS_n)-F_e|\geq|E_{3,4}(BS_n)|-|F_e^0|\geq2(n-2)!-3>0$ for $n\geq5$, $BS_n^{[3,n]}-F_e$ is a subgraph of $H$. Since $|F_e^0|\leq3$, $|V(H)|\geq n!-3$ by Lemma 2.2 (3).

{\bf Subcase 3.2.} {\em $2n-5\leq|F_e^2|\leq4n-13$.}

By Lemma 2.4, $BS_n^2-F_e^2$ has a component $H_2$ with $|V(H_2)|\geq(n-1)!-1$.

Suppose $2n-5\leq|F_e^3|\leq4n-13$. Then $|F_e^0|\leq(8n-21)-(4n-12)-2(2n-5)=1$. Since $|F_e^3|\leq4n-13$, $BS_n^3-F_e^3$ has a component $H_3$ with $|V(H_3)|\geq(n-1)!-1$ by Lemma 2.4. Since $|E_{BS_n}(V(H_i),V(BS_n^4))-F_e|\geq |E_{i,4}(BS_n)|-(|V(BS_n^i)|-|V(H_i)|)-|F_e^0|\geq 2(n-2)!-2-1>0$ for $i\in[1,3]$ and $n\geq5$, $H_i$ is a subgraph of $H$ for every $i\in[1,3]$. If $|V(H_1)|\geq(n-1)!-1$, then $|V(H)|\geq n!-3$. If $|V(H_2)|=(n-1)!$ or $|V(H_3)|=(n-1)!$, then $|V(H)|\geq n!-3$. Now we suppose that $|V(H_1)|=(n-1)!-2$ and $|V(H_2)|=|V(H_3)|=(n-1)!-1$. Let $\{u_{11},u_{12}\}=V(BS_n^1)-V(H_1)$, $u_2\in V(BS_n^2)-V(H_2)$, and $u_3\in V(BS_n^3)-V(H_3)$. Since $|F_e^0|\leq1$, there exists a vertex $v\in(\{u_{11}^+,u_{11}^-,u_{12}^+,u_{12}^-\}-\{u_2,u_3\})\cap V(H)$ such that $(v,u_{11})\in E(BS_n)-F_e$ or $(v,u_{12})\in E(BS_n)-F_e$ by Lemma 2.2 (2). Hence $|V(H)|\geq n!-3$.

Suppose $|F_e^3|\leq 2n-6$. Then $|F_e^0|\leq(8n-21)-(4n-12)-(2n-5)=2n-4$.
By Lemma 2,3, $BS_n^3-F_e^3$ is connected. Since $|E_{3,4}(BS_n)-F_e|\geq|E_{3,4}(BS_n)|-|F_e^0|\geq2(n-2)!-(2n-4)>0$ for $n\geq5$, $BS_n^3-F_e^3$ is a subgraph of $H$. Since $|E_{BS_n}(V(H_i),V(BS_n^4))-F_e|\geq |E_{i,4}(BS_n)|-(|V(BS_n^i)|-|V(H_i)|)-|F_e^0|\geq 2(n-2)!-2-(2n-4)>0$ for $i=1,2$ and $n\geq5$, $H_i$ is a subgraph of $H$. Hence $|V(H)|\geq n!-3$.

{\bf Subcase 3.3.} {\em $|F_e^2|\leq2n-6$.}

By Lemma 2.3, $BS_n^i-F_e^i$ is connected for every $i\in[2,n]$. Since $|E_{i,j}(BS_n)-F_e|\geq|E_{i,j}(BS_n)|-|F_e^0|\geq2(n-2)!-(4n-9)>0$ for $i,j\in[2,n]$ with $i\neq j$ and $n\geq5$, $BS_n^{[2,n]}-F_e$ is a subgraph of $H$. Since $|E_{BS_n}(V(H_1),V(BS_n^{[2,3]}))-F_e|\geq |E_{1,2}(BS_n)|+|E_{1,3}(BS_n)|-2(|V(BS_n^1)|-|V(H_1)|)-|F_e^0|\geq2\times2(n-2)!-2\times2-(4n-9)>0$, $H_1$ is a subgraph of $H$ and $|V(H)|\geq n!-2$.

{\bf Case 4.} {\em $2n-5\leq|F_e^1|\leq4n-13$.}

By Lemma 2.4, $BS_n^1-F_e^1$ has a component $H_1$ with $|V(H_1)|\geq(n-1)!-1$.

{\bf Subcase 4.1.} {\em $|F_e^3|\geq2n-5$.}

In this subcase, $|F_e^0|\leq(8n-21)-3(2n-5)=2n-6$. Since $|E_{i,j}(BS_n)-F_e|\geq|E_{i,j}(BS_n)|-|F_e^0|\geq2(n-2)!-(2n-6)>0$ for $i,j\in[4,n]$ with $i\neq j$ and $n\geq5$, $BS_n^{[4,n]}-F_e$ is a subgraph of $H$.
Since $2n-5\leq|F_e^3|\leq|F_e^2|\leq|F_e^1|\leq4n-13$, $BS_n^i-F_e^i$ has a component $H_i$ with $|V(H_i)|\geq(n-1)!-1$ for $i=2,3$ by Lemma 2.4. Since $|E_{BS_n}(V(H_i),V(BS_n^4))-F_e|\geq |E_{i,4}(BS_n)|-(|V(BS_n^i)|-|V(H_i)|)-|F_e^0|\geq2(n-2)!-1-(2n-6)>0$ for $i\in[1,3]$ and $n\geq5$, $H_i$ is a subgraph of $H$ for every $i\in[1,3]$. Thus $|V(H)|\geq n!-3$.

{\bf Subcase 4.2.} {\em $|F_e^3|\leq2n-6$ and $|F_e^2|\geq2n-5$.}

In this subcase, $|F_e^0|\leq(8n-21)-2(2n-5)=4n-11$. By Lemma 2.3, $BS_n^3-F_e^3$ is connected.
Since $|E_{i,j}(BS_n)-F_e|\geq|E_{i,j}(BS_n)|-|F_e^0|\geq2(n-2)!-(4n-11)>0$ for $i,j\in[3,n]$ with $i\neq j$ and $n\geq5$, $BS_n^{[3,n]}-F_e$ is a subgraph of $H$.
Since $2n-5\leq|F_e^2|\leq|F_e^1|\leq4n-13$, $BS_n^2-F_e^2$ has a component $H_2$ with $|V(H_2)|\geq(n-1)!-1$ by Lemma 2.4. Since $|E_{BS_n}(V(H_i),V(BS_n^4))-F_e|\geq |E_{i,4}(BS_n)|-(|V(BS_n^i)|-|V(H_i)|)-|F_e^0|\geq2(n-2)!-1-(4n-11)>0$ for $i\in[1,2]$ and $n\geq5$, $H_i$ is a subgraph of $H$ for every $i\in[1,2]$. Hence $|V(H)|\geq n!-2$.

{\bf Subcase 4.3.} {\em $|F_e^2|\leq2n-6$.}

In this subcase, $|F_e^0|\leq(8n-21)-(2n-5)=6n-16$.
By Lemma 2.3, both $BS_n^2-F_e^2$ and $BS_n^3-F_e^3$ are connected. We claim $E_{2,3}(BS_n)-F_e\neq\emptyset$ or
$E_{2,4}(BS_n)-F_e\neq\emptyset$; otherwise $|F_e|\geq |E_{2,3}(BS_n)|+|E_{2,4}(BS_n)|=2\times2(n-2)!>8n-21$ for $n\geq5$, a contradiction. Without loss of generality, we assume $E_{2,3}(BS_n)-F_e\neq\emptyset$. Similarly, we can get $E_{2,i}(BS_n)-F_e\neq\emptyset$ or $E_{3,i}(BS_n)-F_e\neq\emptyset$ for every $i\in[4,n]$. Thus $BS_n^{[2,n]}-F_e$ is a subgraph of $H$. Since $|E_{BS_n}(V(H_1),V(BS_n^{[2,3]}))-F_e|\geq |E_{1,2}(BS_n)|+|E_{1,3}(BS_n)|-2(|V(BS_n^1)|-|V(H_1)|)-|F_e^0|\geq2\times2(n-2)!-2\times1-(6n-16)>0$, $H_1$ is a subgraph of $H$. Hence $|V(H)|\geq n!-1$.\q

%{\bf Case 5.} {\em $|F_e^1|\leq2n-6$.}

%In this case, $BS_n^i-F_e^i$ is connected for every $i\in[1,n]$ by Lemma 2.3. We claim that $E_{1,2}(BS_n)-F_e\neq\emptyset$ or $E_{1,3}(BS_n)-F_e\neq\emptyset$; otherwise $|F_e|\geq |E_{1,2}(BS_n)|+|E_{1,3}(BS_n)|=2\times2(n-2)!>8n-21$ for $n\geq5$, a contradiction. Without loss of generality, we assume $E_{1,2}(BS_n)-F_e\neq\emptyset$. Similarly, we can get $E_{1,i}(BS_n)-F_e\neq\emptyset$ or $E_{2,i}(BS_n)-F_e\neq\emptyset$ for every $i\in[3,n]$. Thus $H=BS_n-F_e$ is connected, a contradiction.

\vskip.3cm

\n{\large\bf 3.\quad Edge-fault-tolerant strong Menger edge connectivity of $BS_n$}

\vskip.2cm

We will consider the edge-fault-tolerant strong Menger edge connectivity of $BS_n$ in this section.

\vskip.2cm

{\bf Theorem 3.1} {\em For $n\ge 3$, the bubble-sort star graph $BS_n$ is $(2n-5)$-edge-fault-tolerant strongly Menger edge connected and the bound $2n-5$ is sharp.}

\vskip.2cm
\p Let $F_e\subseteq E(BS_n)$ be an arbitrary faulty edge set with $|F_e|\leq2n-5$. By Lemma 2.3, $BS_n-F_e$ is connected. Let $u,v$ with $u\neq v$ be any two vertices in $BS_n$ and $t=\min\{d_{BS_n-F_e}(u),d_{BS_n-F_e}(v)\}$. By Theorem 1.1, it suffices to show that $u$ and $v$ are connected in $BS_n-F_e-E_f$ for any $E_f\subseteq E(BS_n)-F_e$ with $|E_f|\leq t-1$. Suppose on the contrary, that $u$ and $v$ are disconnected in $BS_n-F_e-E_f$ for some $E_f\subseteq E(BS_n)-F_e$ with $|E_f|\leq t-1$. Since $d_{BS_n-F_e}(u)\leq2n-3$ and $d_{BS_n-F_e}(v)\leq2n-3$, $|E_f|\leq2n-4$. Thus $|F_e\cup E_f|\leq(2n-5)+(2n-4)=4n-9$. By Lemma 2.4, $BS_n-F_e-E_f$ has a component $H$ with $|V(H)|\geq n!-1$. Since $u$ and $v$ are disconnected in $BS_n-F_e-E_f$, $|V(H)|= n!-1$ and $|\{u,v\}\cap V(H)|=1$. Without loss of generality, we assume $u\not\in V(H)$ and $v\in V(H)$. Hence $E_{BS_n}(\{u\},N_{BS_n-F_e}(u))\subseteq E_f$, which implies $|E_f|\geq d_{BS_n-F_e}(u)$, a contradiction to $|E_f|\leq t-1\leq d_{BS_n-F_e}(u)-1$. Hence $BS_n$ is $(2n-5)$-edge-fault-tolerant strongly Menger edge connected.

Next, we will show the bound $2n-5$ is sharp. Let $u,u_1\in V(BS_n)$ with $(u,u_1)\in E(BS_n)$. Let $F_e=E_{BS_n}(u_1)-(u,u_1)$ and $v\in V(BS_n)-N_{BS_n}(u_1)-\{u_1\}$ (see Fig.2). Then $|F_e|=2n-4$, $d_{BS_n-F_e}(u)=d_{BS_n-F_e}(v)=2n-3$. Obviously, there are at most $2n-4$ edge-disjoint $(u,v)$-paths.\q
\begin{figure}[ht]
   \begin{center}
    \includegraphics[scale=0.6]{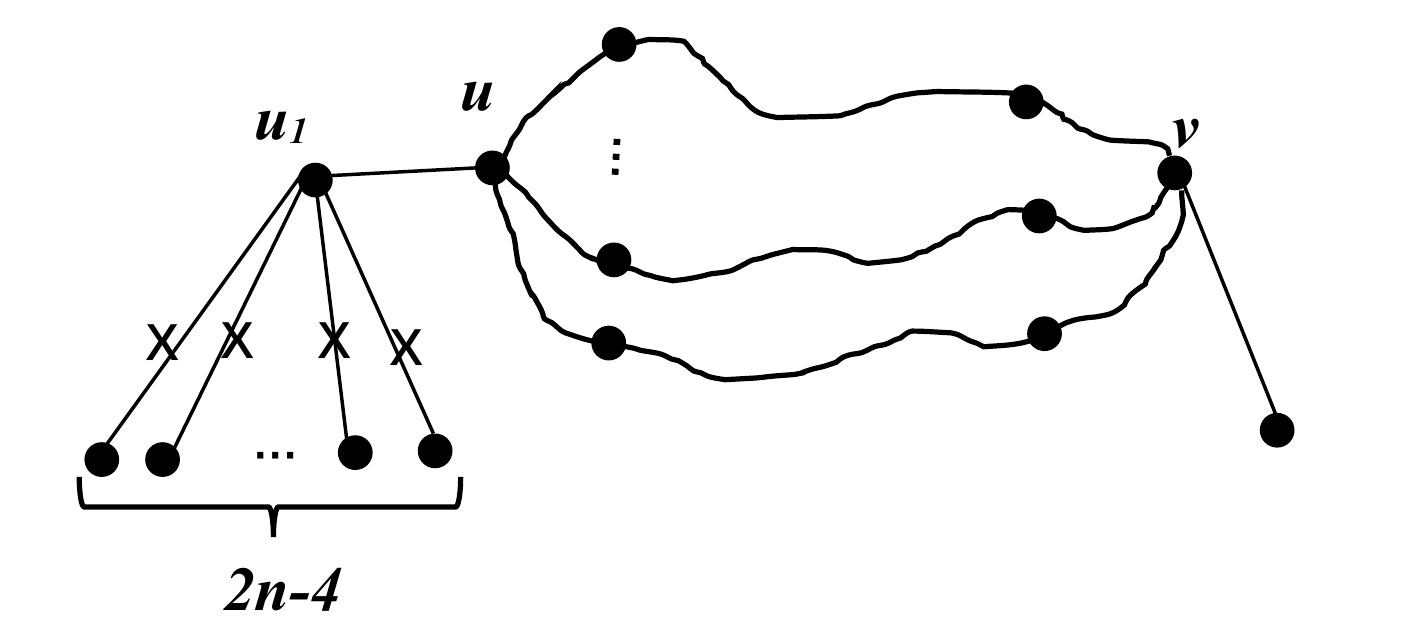}
    \end{center}
    \vspace{-0.5cm}\caption{\label{F2-5} Illustration of Theorem 3.1.}
\end{figure}

\n{\large\bf 4.\quad Conditional edge-fault-tolerant strong Menger edge connectivity of $BS_n$}
\vskip.2cm

We will consider the conditional edge-fault-tolerant strong Menger edge connectivity of $BS_n$ in this section.

\vskip.2cm

{\bf Theorem 4.1} {\em For $n\ge 4$, the bubble-sort star graph $BS_n$ is $(6n-17)$-conditional edge-fault-tolerant strongly Menger edge connected and the bound $6n-17$ is sharp.}

\vskip.2cm
\p Let $F_e\subseteq E(BS_n)$ be an arbitrary faulty edge set with $|F_e|\leq6n-17$ and $\delta(BS_n-F_e)\geq2$. Since $|F_e|\leq6n-17\leq6n-14$ and $\delta(BS_n-F_e)\geq2$, $BS_n-F_e$ is connected by Lemma 2.7. Let $u,v$ with $u\neq v$ be any two vertices in $BS_n$ and $t=\min\{d_{BS_n-F_e}(u),d_{BS_n-F_e}(v)\}$. By Theorem 1.1, it suffices to show that $u$ and $v$ are connected in $BS_n-F_e-E_f$ for any $E_f\subseteq E(BS_n)-F_e$ with $|E_f|\leq t-1$. Suppose on the contrary, that $u$ and $v$ are disconnected in $BS_n-F_e-E_f$ for some $E_f\subseteq E(BS_n)-F_e$ with $|E_f|\leq t-1$. Since $d_{BS_n-F_e}(u)\leq2n-3$ and $d_{BS_n-F_e}(v)\leq2n-3$, $|E_f|\leq2n-4$. Thus $|F_e\cup E_f|\leq(6n-17)+(2n-4)=8n-21$. By Lemma 2.9, $BS_n-F_e-E_f$ has a component $H$ with $|V(H)|\geq n!-3$. Since $u$ and $v$ are disconnected in $BS_n-F_e-E_f$, $|\{u,v\}\cap V(H)|\leq1$. Without loss of generality, we assume $u\not\in V(H)$. Let $H_1$ be the component in $BS_n-F_e-E_f$ containing $u$. If $d_{H_1}(u)=0$, then $E_{BS_n}(\{u\},N_{BS_n-F_e}(u))\subseteq E_f$, which implies $|E_f|\geq d_{BS_n-F_e}(u)$, a contradiction to $|E_f|\leq t-1\leq d_{BS_n-F_e}(u)-1$. Suppose that $d_{H_1}(u)=i$ ($i\in[1,2]$). Since $BS_n$ is bipartite, $H_1$ is a path $P_2$ or $P_3$ and there are $i$ vertices in $V(H_1)-\{u\}$ that have degree one in $H_1$. Since $\delta(BS_n-F_e)\geq2$, every vertex with degree one in $H_1$ is incident with at least one edge in $E_f$. Thus $|E_f|\geq d_{BS_n-F_e}(u)-i+i=d_{BS_n-F_e}(u)$, a contradiction to $|E_f|\leq t-1\leq d_{BS_n-F_e}(u)-1$. Hence $BS_n$ is $(6n-17)$-conditional edge-fault-tolerant strongly Menger edge connected.

Next, we will show the bound $6n-17$ is sharp. Let $u,u_1,u_2,u_3\in V(BS_n)$ with $(u,u_1),(u_1,u_2),(u_2,u_3),(u_3,u)\in E(BS_n)$ and $u_{11}\in N_{BS_n}(u_1)-\{u,u_2\}$. Let $F_e=E_{BS_n}(u_1)\cup E_{BS_n}(u_2)\cup E_{BS_n}(u_3) -\{(u,u_1),(u_1,u_2),(u_2,u_3),(u_3,u),(u_1,u_{11})\}$ and $v\in V(BS_n)-N_{BS_n}(u_1)\cup N_{BS_n}(u_2)\cup N_{BS_n}(u_3)$ (see Fig.3). Then $|F_e|=(2n-6)+2(2n-5)=6n-16$, $d_{BS_n-F_e}(u)=d_{BS_n-F_e}(v)=2n-3$, and $\delta(BS_n-F_e)\geq2$ for $n\geq4$. Obviously, there are at most $2n-4$ edge-disjoint $(u,v)$-paths.\q
\begin{figure}[ht]
   \begin{center}
    \includegraphics[scale=0.6]{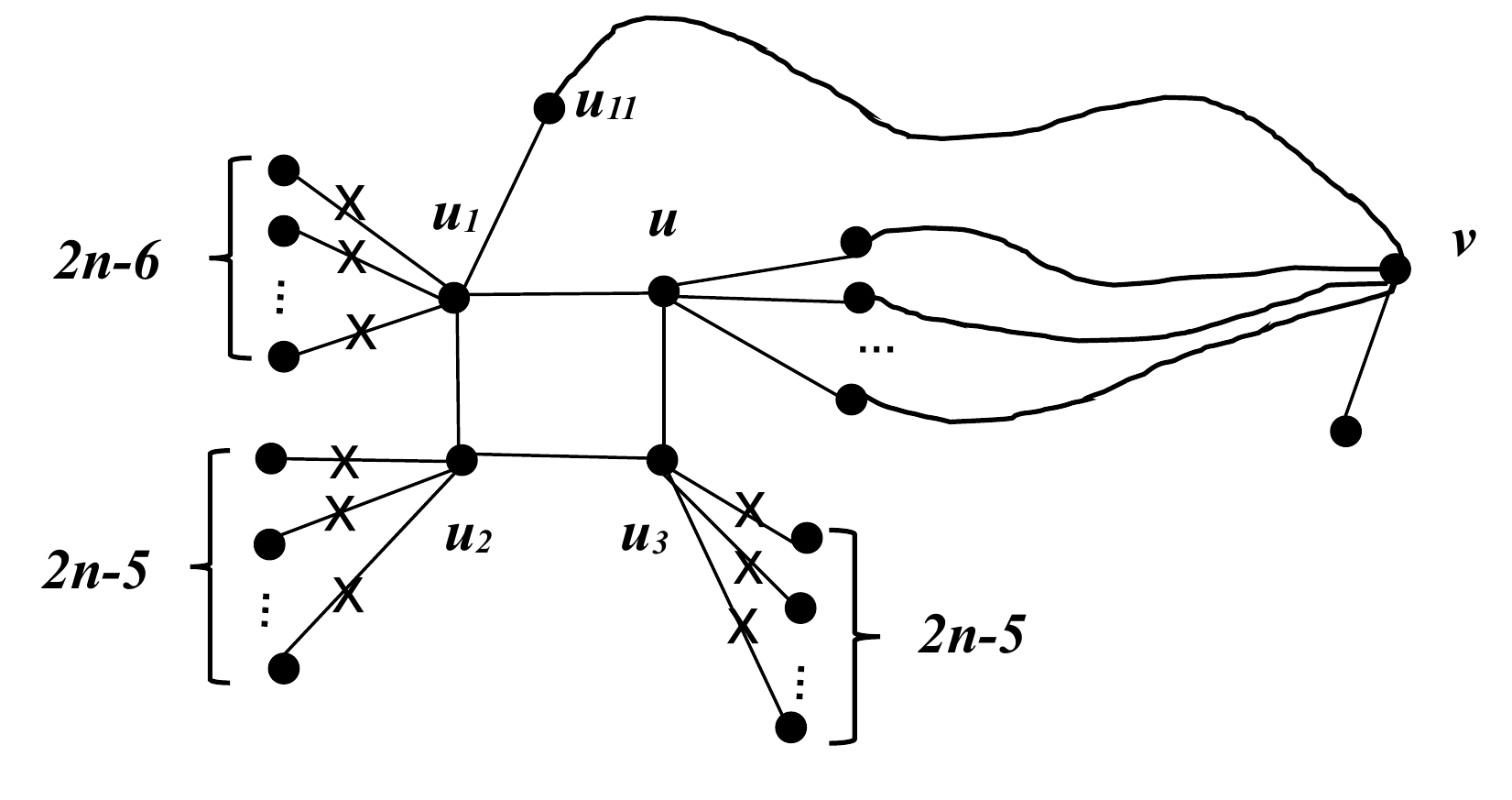}
    \end{center}
    \vspace{-0.5cm}\caption{\label{F2-5} Illustration of Theorem 4.1.}
\end{figure}

\n{\large\bf 5.\quad Conclusion}

\vskip.2cm
In this paper, we study the edge-fault-tolerant strong Menger edge connectivity of $n$-dimensional bubble-sort star graph $BS_n$. We show that every pair of distinct vertices $u$ and $v$ in $BS_n$ are connected by $\min\{d_{BS_n-F_e}(u),d_{BS_n-F_e}(v)\}$ edge-disjoint paths in $BS_n-F_e$, where $F_e$ is an arbitrary edge subset of $BS_n$ with $|F_e|\leq 2n-5$. We also show that every pair of distinct vertices $u$ and $v$ in $BS_n$ are connected by $\min\{d_{BS_n-F_e}(u),d_{BS_n-F_e}(v)\}$ edge-disjoint paths in $BS_n-F_e$, where $F_e$ is an arbitrary edge subset of $BS_n$ with $|F_e|\leq 6n-17$ and $\delta(BS_n-F_e)\geq2$.
Moreover, we give two examples to show that our results are optimal. The connectivity and edge connectivity of interconnection network determine the fault tolerance of the network. They are issues worth studying.

\vskip.3cm

\n{\large\bf Acknowledgements}
\vskip.2cm
This research is supported by National Natural Science Foundation of China (No. 11801450), Natural Science Foundation of  Shaanxi Province, China (No. 2019JQ-506).

\end{document}